\documentclass{amsart}	
\usepackage{amsmath, amsthm, amssymb, amscd}
\usepackage[all]{xy}

\newcommand{\Z}{\mathbb{Z}}

\newcommand{\bbC}{\mathbb{C}}

\newcommand{\ZZ}{\mathbb{Z}}
\newcommand{\CC}{\mathbb{C}}

\newcommand{\A}{\mathcal{A}}
\newcommand{\cA}{\mathcal{A}}

\newcommand{\K}{\mathcal{K}}
\newcommand{\cL}{\mathcal{L}}
\newcommand{\cC}{\mathcal{C}}

\theoremstyle{plain}
\newtheorem{theorem}{Theorem}[section]

\newtheorem{corollary}[theorem]{Corollary}

\theoremstyle{remark}

\theoremstyle{definition}

\begin{document}
\date{}
\title{Entire cyclic homology of stable continuous trace algebras}

\author[V. Mathai]{Varghese Mathai}
\address[Varghese Mathai]
{Department of Pure Mathematics\\
University of Adelaide\\
Adelaide, SA 5005 \\
Australia}
\email{vmathai@maths.adelaide.edu.au }

\author[D. Stevenson]{Danny Stevenson}
\address[D. Stevenson]{Department of Mathematics\\
University of California\\
Riverside, CA, USA}
\email{dstevens@maths.adelaide.edu.au }

\thanks{The authors acknowledge the support of the Australian
Research Council.}


\begin{abstract}
A central result here is the computation of the
entire cyclic homology of canonical smooth subalgebras
of stable continuous trace $C^*$-algebras
having smooth manifolds $M$ as their spectrum.
More precisely, the entire cyclic homology is shown to be canonically
isomorphic to the continuous periodic cyclic homology for these
algebras. By an earlier result of the authors, one concludes
that the entire cyclic homology of the algebra is canonically
isomorphic to the twisted de Rham cohomology of $M$.
\end{abstract}


\maketitle


\section{Introduction}

Entire cyclic homology was introduced by Connes in \cite{Connes2}, as a
version of cyclic homology that is better suited to study the dual space
of higher rank discrete groups
and also certain infinite dimensional spaces occurring in the study of
constructive field theory. Entire cyclic homology was placed in a universal
setting by Meyer in \cite{Meyer}, which is what we use here.

Recently, Brodzki and Plymen \cite{BP} computed the entire cyclic
homology of the Schatten class ideals $\cL^p$, for $p\ge 1$. More
precisely, they proved that the entire cyclic
homology, $HE_\bullet(\cL^p)$ is isomorphic to the continuous periodic
cyclic homology $HP_\bullet^{\mathrm{cont}}(\cL^p)$, which had been
previously computed by
Cuntz \cite{Cu1, Cu3} to be equal to  $HP_\bullet^{\mathrm{cont}}(\bbC)$. 
On the other hand, Puschnigg \cite{Puschnigg}
had computed the entire cyclic homology of the Fr\'echet algebra
of smooth functions on a compact manifold, by showing using
excision that it was isomorphic to the continuous periodic cyclic homology,
which had been computed earlier by Connes in \cite{Connes}.

In this paper, we use excision in entire cyclic homology
\cite{Puschnigg, Meyer}
and in continuous periodic cyclic homology \cite{Cu2, CuQu},
    together with the computation
of the entire cyclic homology of  the Schatten class ideals \cite{BP}, to
reduce the computation of the entire cyclic homology of
a canonical smooth subalgebra
of stable continuous trace $C^*$-algebras
having smooth manifolds $M$ as their spectrum, to the computation of
the continuous periodic cyclic homology of these algebras.
By an earlier result of the authors \cite{MS}, one concludes
that the entire cyclic homology of the algebra is canonically
isomorphic to the twisted de Rham cohomology of $M$.

\section{Preliminaries}

A collection $\mathfrak{B}$ of subsets of a vector space $V$ is called a
\emph{bornology} if,
roughly speaking, it is
closed under forming finite sums and taking subsets.  For more
details we refer to \cite{Meyer}.  A vector space $V$ is called
a \emph{bornological vector space} if it comes equipped with a
bornology $\mathfrak{B}$.  We will usually denote a bornological
vector space by $V$ if the bornology $\mathfrak{B}$ is understood;
if any confusion could occur we will write $(V,\mathfrak{B})$.
A subset $B\subset V$ of a bornological vector space $V$ is called
\emph{bounded} if $B\in \mathfrak{B}$.
Recall that for any bornological vector
space $V$, the operations of addition $V\times V \to V$
and scalar multiplication $\CC\times V\to V$ are compatible with the
bornological structure of $V$ in the sense that these are bounded linear
maps, where a bounded map in this context means a map which
sends bounded sets
to bounded sets in the respective bornologies.  Given a
vector subspace $W\subset V$ of a bornological vector
space $V$ we can define the \emph{subspace} bornology
on $W$ and also the \emph{quotient} bornology on $V/W$.
For more details we refer to \cite{Meyer}.
A \emph{bornological algebra}
is an algebra $A$ equipped with a bornology making it into a
bornological
vector space and is such that the product $A\times A\to A$ is
bounded.
Given a bornological vector space $V$, we can form its \emph{completion}
$V^c$, this is characterised by the usual universal property: any
bounded
linear map $V\to W$ to a complete bornological vector space $W$
factorises
uniquely through $V^c$ via a bounded linear map $V\to V^c$.  Given
two bornological vector spaces $V$ and $W$ we can equip the algebraic
tensor product $V\otimes W$ with the \emph{tensor product} bornology.
This is the convex bornology on $V\otimes W$ generated by sets of the
form $B_1\otimes B_2$ where $B_1\in \mathfrak{B}_V$ and
$B_2\in \mathfrak{B}_W$.
We denote the
completion of $V\otimes W$ in the tensor product bornology
as $(V\hat{\otimes} W)^c$ or simply $V\hat{\otimes}W$.  If $A$ and
$B$ are bornological algebras then the completed bornological
tensor product $A\hat{\otimes}B$ is also a bornological algebra.
If $V$ is a Fr\'{e}chet space then there is associated to $V$ a
canonical
bornology called the \emph{pre-compact} bornology which makes
$V$ into a complete bornological vector space.  For this
canonical bornology we take the collection
of all pre-compact subsets of $V$.  If $A$ is a Fr\'{e}chet algebra
then the
pre-compact bornology makes $A$ into a complete bornological algebra.
Continuous maps between Fr\'{e}chet spaces exactly correspond to
bounded linear maps on the associated bornological spaces.  If
$W\subset V$
is a closed subspace of the Fr\'{e}chet space $V$ then the pre-compact
bornology
on the quotient $V/W$ coincides with the quotient bornology on $V/W$.
One final property of Fr\'{e}chet spaces we will need is that the
\emph{bornological
tensor product} $V\hat{\otimes}W$ of two Fr\'{e}chet spaces equipped
with the pre-compact bornology coincides with the pre-compact
bornology on the projective tensor
product $V\hat{\otimes}_{\pi}W$ of $V$ and $W$.

We now recall the definition of entire and periodic cyclic
homology following \cite{Meyer, Cu1, CuQu, Perrot}.  Let
$A$ be a complete bornological algebra.  Recall
that one first forms the $\ZZ_2$-graded vector space
$$
\Omega A = \bigoplus_{n=1}^\infty \Omega^n A
$$
where $\Omega^n A = \tilde{A}\hat{\otimes}A^{\hat{\otimes}\, n}$ for
$n\geq 1$.
with $\tilde{A} = A\oplus \CC$ the unitisation of $A$.  We set
$\Omega^0 A = A$.  $\Omega A$ is
made into a bornological algebra using the \emph{entire} bornology
\cite{Meyer, Perrot}, i.e. the bornology generated by sets of the
form $\cup_{n\geq 0}[n/2]!\tilde{S}dS^n$, where $S\in \mathfrak{B}$.
The completion of $\Omega A$ in this bornology is denoted by
$\Omega_{\epsilon} A$.  It can be shown that the standard operators
$b$ and $B$ are bounded and extend to $\Omega_{\epsilon} A$.  The
\emph{entire cyclic homology} $HE_*(A)$ of $A$ is defined to be the
homology of the complex $(\Omega_{\epsilon} A, b + B)$.  That this
definition of $HE_*(A)$ is the same as the definition given in
\cite{Meyer}
using the $X$-complex of Cuntz and Quillen is explained in
\cite{Perrot}.

We can also form the periodic complex $\widehat\Omega A =
\displaystyle{\prod_{n=0}^\infty}
\Omega^n A$ equipped with the direct product bornology.  Again the
operators $b$ and $B$ are bounded and so we can form the homology
of the complex $(\widehat\Omega A,b+B)$.  This homology is the (bornological)
periodic cyclic homology of $A$.  Notice that the
canonical map $\Omega_{\epsilon}A\to \widehat\Omega A$ is bounded, thus we have
a
morphism of complexes and hence a canonical map $HE_*(A)\to HP_*(A)$.
As
pointed out in \cite{Meyer}, if $A$ is a Fr\'{e}chet algebra equipped
with the
pre-compact bornology, then the bornological periodic cyclic homology
and the continuous periodic cyclic homology
$HP^{\mathrm{cont}}_*(A)$ of $A$ coincide.  Therefore
we have a natural map $HE_*(A)\to HP^{\mathrm{cont}}_*(A)$.
It can be shown that this map is compatible with the connecting
homomorphisms in the long exact sequences in the entire and
periodic cyclic theories associated to extensions $0\to A\to B\to C\to  
0$
of bornological algebras equipped with a continuous linear splitting.

\section{Main Theorems}

The following is a generalization of Theorem 6.1 in \cite{Puschnigg},
which is possible thanks to the main result in \cite{BP}.

\begin{theorem}\label{thm:class}
Let $\cC$ be the smallest class of Fr\'echet algebras satisfying the  
following properties:
\begin{enumerate}
\item $\bbC\in \cC$, $\cL^p \in \cC$ for all $p\ge 1$;
\item $\cC$ is closed under smooth homotopy equivalence;
\item If in an extension admitting a continuous linear section, two of  
the
algebras
belong to $\cC$, then so does the third.
\end{enumerate}

Then for every algebra $\cA \in \cC$, the natural map between
the
continuous periodic homology and the entire cyclic homology induces
isomorphisms
$$
HE_\bullet (\cA) \cong  HP^{\mathrm{cont}}_\bullet (\cA)
$$
where $\cA$ is equipped with the precompact bornology.
\end{theorem}

\begin{proof} Consider the class $\cC'$
of Fr\'echet algebras endowed with the canonical
complete bornology,
such that transformations between continuous periodic
homology and entire cyclic homology are isomorphisms.
Then by \cite{Connes2}, $\cC'$ contains $\bbC$, and by
\cite{BP}, $\cC'$ contains $\cL^p$ for $p\ge 1$.

Since continuous periodic homology and
entire cyclic homology are smooth homotopy functors see \cite{Connes, Meyer}, the
class $\cC'$ is closed under smooth homotopy equivalence.
The forgetfulness map from entire cyclic homology to continuous periodic
cyclic homology, $HE_\bullet \to HP^{\mathrm{cont}}_\bullet$ is compatible with the
long exact sequences associated to extensions of Fr\'echet algebras,
as was observed in the previous section.
So, if in an extension admitting a continuous linear section, two of the
algebras
belong to $\cC'$, then by the excision theorem, \cite{Meyer}
together with the five lemma imply that the third algebra in the
extension
also belongs to $\cC'$.
Thus $\cC \subset \cC'$, proving the theorem.
\end{proof}

Let $A$ be a stable continuous trace $C^*$-algebra with spectrum
a smooth, compact manifold $M$. By a fundamental theorem of
Dixmier-Douady,
we know that $A =C(M, \K(P) )$ is the algebra of continuous sections of
a smooth,
locally trivial bundle $\K(P) = P \times_{PU} \K$ on $M$ with fibre the
algebra $\K$ of compact operators on a separable
Hilbert space associated to a principal $PU$
bundle $P$ on $M$ via the adjoint action of
$PU$ on $\K$.  Here $PU$ is the group
of projective unitary operators on the Hilbert space.
Such algebras  $A$
are classified up to isomorphism by their
Dixmier-Douady invariant   $\delta(P)\in H^3(M;\Z)$.
Inside $A$ we can consider
a dense, canonical smooth $*$-subalgebras $\A_p = C^\infty(M,
\cL^p(P))$,
$p \ge 1$, consisting of all
smooth sections of the sub-bundle $\cL^p(P) = P \times_{PU} \cL^p$ of
$\K(P)$ with
fibre the Schatten class ideal $\cL^p$ of operators on the Hilbert space
and structure group $PU$.

\begin{theorem}\label{thm:main}
Let $M$ be a smooth connected compact manifold and $P\to M$ a principal
$PU$ bundle. Consider the nuclear Fr\'echet-algebra of smooth sections
$\cA_p$, $p\ge 1$, of the bundle
    $\cL^p(P) = P \times_{PU} \cL^p$. Then for all $p\ge 1$,
    there are natural isomorphisms,
\begin{equation}
HE_\bullet (\cA_p) \cong HP^{\mathrm{cont}}_\bullet
(\cA_p).
\end{equation}
\end{theorem}

\begin{proof}
The proof that we present is along the lines of Theorem 6.2 in
\cite{Puschnigg}
and the difference is only in the details.

Let $N$ be a smooth compact submanifold of $M$ of codimension one,
with trivial normal bundle. Let $C^\infty(M, N, \cL^p(P))$ be the
algebra of smooth sections of $\cL^p(P)$ vanishing when restricted to
$N$ and let $C^\infty_0(M, N, \cL^p(P))$ be the algebra of smooth
sections of $\cL^p(P)$ vanishing to infinite order when restricted to
$N$. Then the inclusion map
\begin{equation}
\iota : C^\infty_0(M, N, \cL^p(P)) \hookrightarrow C^\infty(M, N,
\cL^p(P))
\end{equation}
is a smooth homotopy equivalence. To see this, let $\phi \in {\rm
Diff}({\mathbb R})$
be a diffeomorphism which is equal to the identity outside of the
interval $[-1, 1]$ and has the
property that it vanishes to infinite order at the origin. Any open
tubular neighbourhood
$W$ of $N$ in $M$ can be identified with $N\times {\mathbb R}$, since
$N$ is assumed to have
trivial normal bundle in $M$. Now we can extend the diffeomorphism
${\rm Id} \otimes \phi$ of
$N\times {\mathbb R}$, to a diffeomorphism $\Phi$ of $M$,
by setting it equal to the identity outside of $W$. The algebra
homomorphism $\Phi^*: C^\infty(M,  \cL^p(P)) \to C^\infty(M,
\Phi^*\cL^p(P))$ has the property that
it maps $C^\infty(M, N, \cL^p(P))$ to $C^\infty_0(M, N,
\Phi^*\cL^p(P))$. Since ${\rm Id} \otimes \phi$
is homotopic to the identity, it follows that $\Phi$ is also homotopic to
the identity, so that the bundles
$ \Phi^*\cL^p(P)$ and $\cL^p(P)$ are isomorphic. Finally, it is clear
that $\Phi^*$ is an
inverse to $\iota$ up to smooth homotopy.

Next we prove by induction over the dimension that $C^\infty(M,
\cL^p(P))  \in \cC$,
where the class $\cC$ is as defined in Theorem \ref{thm:class}. In
dimension zero,
$C^\infty(M,  \cL^p(P)) =  \cL^p$ and
the main result in \cite{BP} establishes the claim in this case. So we
will now assume that
${\rm dim}(M)>0$. Let $f: M \to \mathbb R^+$ be a Morse function with
finitely many
critical points $x_1, x_2, \ldots x_k$ lying on pairwise different
level surfaces of $f$.
Set $N_i = f^{-1}(t_i)$, $i=1, \ldots k$ be level surfaces of $f$
having the property that
$t_{i-1} < f(x_i) < t_i$, $i=1, \ldots k$. Then the extension of
nuclear Fr\'echet-algebras
$$
0\to C^\infty(M,  \bigsqcup N_i, \cL^p(P)) \to C^\infty(M,  \cL^p(P))
\to C^\infty(\bigsqcup N_i, \cL^p(P)) \to 0
$$
possesses a bounded linear section, so by the induction hypothesis, it
suffices
to verify that $C^\infty(M,  \bigsqcup N_i, \cL^p(P)) \in \cC$. As seen
earlier, there are
homotopy equivalences {\small
$$
\begin{array}{lcl}
C^\infty(M,  \bigsqcup N_i, \cL^p(P)) &\cong& C^\infty_0(M,  \bigsqcup
N_i, \cL^p(P))\\[+7pt]
&\cong & \bigoplus C^\infty_0(f^{-1}([t_i, t_{i+1}]), N_i\bigsqcup
N_{i+1},  \cL^p(P)|_{f^{-1}([t_i, t_{i+1}])})
\end{array}
$$}
So it suffices to verify that each summand on the right hand side is in
$\cC$.

Now $ C^\infty_0(f^{-1}([t_i, t_{i+1}]), N_i\bigsqcup N_{i+1},
\cL^p(P)|_{f^{-1}([t_i, t_{i+1}])})$
is smoothly homotopy equivalent to
$ C^\infty_0(D^{n_i} \times S^{m_i}, \partial D^{n_i} \times S^{m_i},
\cL^p(P)|_{D^{n_i} \times S^{m_i}})$
which is smoothly homotopy equivalent to
$ C^\infty (D^{n_i} \times S^{m_i}, \partial D^{n_i} \times S^{m_i},
\cL^p(P)|_{D^{n_i} \times S^{m_i}})$
where $n_i + m_i = n$ and $n_i>0$, for $i=1, \ldots , k$. Now the
extension of nuclear Fr\'echet-algebras
\begin{multline}
0\to  C^\infty (D^{n_i} \times S^{m_i}, \partial D^{n_i} \times
S^{m_i},  \cL^p(P)|_{D^{n_i} \times S^{m_i}})
\to  \\ C^\infty (D^{n_i} \times S^{m_i},  \cL^p(P)|_{D^{n_i} \times
S^{m_i}})
\to  C^\infty ( \partial D^{n_i} \times S^{m_i},  \cL^p(P)|_{\partial
D^{n_i} \times S^{m_i}})\to 0
\end{multline}
together with the smooth homotopy equivalence  of $C^\infty (D^{n_i}
\times S^{m_i},  \cL^p(P)|_{D^{n_i} \times S^{m_i}})$ and $C^\infty
(S^{m_i},  \cL^p(P)|_{S^{m_i}})$ and the induction hypothesis
establishes that $$ C^\infty_0(f^{-1}([t_i, t_{i+1}]), N_i\bigsqcup
N_{i+1},  \cL^p(P)|_{f^{-1}([t_i, t_{i+1}])})
\in \cC,$$for all $i=1, \ldots , k$. The induction is complete and
therefore $\A_p = C^\infty(M, \cL^p(P))
\in \cC$.

\end{proof}

\begin{corollary}
The entire cyclic homology $HE_\bullet(\A_p)$ is isomorphic to the
twisted de Rham cohomology $H^\bullet(M;c(P))$
for some closed $3$-form $c(P)$ on $M$ such that
$\frac{1}{2\pi i}c(P)$ represents the image of
the Dixmier-Douady invariant $\delta(P)$
in real cohomology.
\end{corollary}

\begin{proof}
By Theorem \ref{thm:main}, we know that the entire cyclic
homology $HE_\bullet(\A_p)$ is isomorphic to the continuous
periodic cyclic homology $HP_\bullet(\A_p)$. By the main
result in \cite{MS}, the continuous
periodic cyclic homology $HP_\bullet(\A_p)$ is isomorphic to the
twisted de Rham cohomology $H^\bullet(M;c(P))$
for some closed $3$-form $c(P)$ on $M$ such that
$\frac{1}{2\pi i}c(P)$ represents the image of
the Dixmier-Douady invariant $\delta(P)$
in real cohomology.
\end{proof}


\end{document}